\documentclass[reqno,12pt]{amsart} 
\usepackage{amsmath, amsfonts, amsthm, amssymb, amscd}
\usepackage[mathcal]{euscript} 
\usepackage{dsfont, pxfonts, mathrsfs, accents, verbatim} 
\theoremstyle{plain}
        \newtheorem{theorem}{Theorem}[section]
	         
        \newtheorem{lemma}[theorem]{Lemma}

        \newtheorem{remark}[theorem]{Remark}  
\numberwithin{equation}{section} 
\newcommand \RR 		{\mathbb{R}}  
  
\newcommand \del 		\partial
\newcommand \eps 		\epsilon 
\newcommand \be 		{\begin{equation}}
\newcommand \ee 		{\end{equation}} 
\newcommand \Hcal 	 	{\mathcal H}
 
\newcommand \hatl 		{\widehat l}
\newcommand \hatr 		{\widehat r}
\newcommand \lam		{\lambda} 
\newcommand \hatlam		{\widehat\lambda} 
\newcommand \la 		\langle 
\newcommand \ra 		\rangle 
\newcommand \diag 		{\text{diag}} 
\begin{document}
\title[Phase dynamics with physical viscosity and capillarity]
{Singular limits in phase dynamics
\\
with physical viscosity and capillarity
}
\author
%   Note for the Galley Proofs : the name "LeFloch" should appear as a single word with the third letter capitalized. 
   [K.T. Joseph and P.G. LeFloch]
    {K.T. Joseph and P.G. LeFloch}
\address
   {K.T. Joseph\\
   School of Mathematics, 
Tata Institute of Fundamental Research, Homi Bhabha Road,
Mumbai 400005, India. E-mail: {ktj@math.tifr.res.in}
\\
\newline 
P.~G. LeFloch\\
   Laboratoire J.-L. Lions \& CNRS UMR 7598\\
   University of Paris VI, 75252 Paris, France. E-mail : {lefloch@ann.jussieu.fr}  
}
\date{August 26, 2006.}



\begin{abstract} 
Following pioneering work by Fan and Slemrod who studied the effect of artificial viscosity terms, 
we consider the system of conservation laws arising in liquid-vapor phase dynamics
with {\sl physical} viscosity and capillarity effects taken into account. 
Following Dafermos we consider self-similar solutions to the Riemann problem and 
establish uniform total variation bounds, allowing us to deduce new existence results.  
Our analysis cover both the hyperbolic and the hyperbolic-elliptic regimes and apply to arbitrarily large 
Riemann data. 

The proofs rely on a new technique of reduction to two coupled scalar equations
associated with the two wave fans of the system.  Strong $L^1$ convergence to a weak solution of bounded variation 
is established in the hyperbolic regime, while in the hyperbolic-elliptic regime
a stationary singularity near the axis separating the two wave fans, or more generally  
an almost-stationary oscillating wave pattern (of thickness depending upon the capillarity-viscosity ratio)
are observed which prevent the solution to have globally bounded variation.  
\end{abstract}

{\sl \small \noindent To appear in : Proc. A Royal Soc. Edinburgh.}

\

\maketitle  

%\tableofcontents

%=============================================================== 

\section{Introduction}

The Navier-Stokes equations for van der Waals fluids with viscosity and capillarity effects included, 
allow one to model the dynamics of liquid-vapor flows. The associated set of first-order conservation laws 
is of hyperbolic or hyperbolic-elliptic type, and admits propagating discontinuities (shock waves). 
The singular limit corresponding to vanishing viscosity and capillarity coefficients allows one to select physically
admissible, discontinuous solutions to the first-order conservation laws. In particular, this yields an approach to 
select solutions to the so-called Riemann problem, when the initial data consists of two constant states
separated by a single jump.  Our purpose in the present paper is to derive uniform (viscosity-capillarity independent)
bounds and to justify this singular limit in the context of nonlinear elasticity and phase transition dynamics. 
 
Attention will be concentrated on the following system of two conservation laws: 
\be
\label{elastic}
\aligned 
& {\del v \over \del t} - {\del \over \del x} \Big( \sigma(w) - \eps \, {\del v \over \del x} + \delta {\del^2 w \over \del x^2} \Big) = 0,   
\\
& {\del w \over \del t} - {\del v \over \del x} = 0,
\endaligned 
\ee
where $v=v(t,x)$ and $w=w(t,x)$ represent the velocity and deformation gradient 
of the fluid or solid material under consideration. 
We consider also the associated first-order system
\be
\label{elastic-first}
\aligned 
& {\del v \over \del t} - {\del \over \del x}( \sigma(w) ) = 0,   
\\
& {\del w \over \del t} - {\del v \over \del x} = 0,
\endaligned 
\ee
and impose Riemann initial data 
\be
\label{boundary}
(v,w)(0,x) = \begin{cases}
v_l, w_l,    & x < 0, 
\\
v_r, w_r,    & x > 0, 
\end{cases}
\ee
where $v_l, v_r, w_l, w_r \in \RR$ are given constants. Solutions of \eqref{elastic-first}-\eqref{boundary} 
are known to be self-similar, that is, to depend only upon the variable $y:=x/t$. 

When viscosity and capillarity terms are taken into account, 
the corresponding set of differential equations reads 
\be
\label{elastic2}
\aligned 
& -y \, v' - \sigma(w)' = \eps \, v'' - \delta w''', 
\\
& - y \, w' - v' = 0,
\endaligned 
\ee
supplemented with the boundary conditions 
\be
\label{boundary1}
\lim_{y \to -\infty} (v,w)(y) = (v_l, w_l), 
\qquad 
\lim_{y \to +\infty} (v,w)(y) = (v_r, w_r).  
\ee
Our main objective in this paper is to establish : 
\begin{enumerate}
\item the existence of a smooth, self-similar solution to \eqref{elastic2}-\eqref{boundary1} having uniformly bounded total variation,  
\be
\label{BVbound}
TV(v_\eps, w_\eps) \lesssim |v_r - v_l| + |w_r - w_l|, 
\ee  
with the implied constant being {\sl independent} of $\eps, \delta$ within the range $\delta /\eps^2 <<1$, and 
\item the strong convergence of $v_\eps, w_\eps$ toward a weak, self-similar 
solution $v,w$ to the first-order conservation 
laws \eqref{elastic-first}. 
\end{enumerate} 

An outline of our main results is as follows. 

Assuming first that the system \eqref{elastic-first} is {\sl hyperbolic} we provide a rather direct proof of the 
above two properties, first in the viscosity-only case (Section~2) and then for general viscosity-capillarity (Section~3). 
In this regime, it is known that the limiting solutions in general contains nonclassical shocks, with 
depend upon the ratio $\delta/\eps^2$; see \cite{LeFloch2}. 

Next, we investigate the generalization to the hyperbolic-elliptic regime, and discover a concentration phenomena 
near the axis $x=0$ separating the two wave fans. This seems to be consistent with numerical experiments 
with the model under consideration, but it would be interesting to check numerically the feature discovered here analytically. 
In the context of phase dynamics, it is also well-known that the limiting solutions contains subsonic 
phase boundaries which again depend on the capillarity to viscosity ratio; 
see \cite{Slemrod1,Slemrod2,Truskinovsky,AbeyaratneKnowles,LeFloch1,Shearer,ShearerYang}. 

Observe that our results cover arbitrary large Riemann data and physical viscosity and capillarity terms.  
Our results supplement the earlier, pioneering 
work by Fan and Slemrod  \cite{FanSlemrod,FanSlemrod2} where an artificial, ``full'' viscosity was used. 
Our technique of proof in the present paper is quite different from the one in \cite{FanSlemrod}, 
as we introduce a decomposition of the system of equations \eqref{elastic2} into two coupled scalar equations. 

Finally, we provide various generalization to the boundary value problem
and more general classes of conservation laws. 

Recall that the activity on self-similar vanishing viscosity limits started with an extensive research by 
Dafermos \cite{Dafermos1}--\cite{Dafermos4} (see also \cite{DafermosDiPerna})
who advocated the use of self-similar regularizations to capture the whole structure 
of wave fans within solutions of the Riemann problem. Self-similar approximations in the context of phase transition dynamics
were studied by Slemrod \cite{Slemrod1,Slemrod2}, Fan and Slemrod \cite{FanSlemrod}, 
who covered large data solution and artificial regularization terms.
Next, small data solutions to general systems 
were treated by Tzavaras \cite{Tzavaras} (conservative systems) and LeFloch and Tzavaras 
\cite{LeFlochTzavaras1,LeFlochTzavaras2} (nonconservative systems). 

Self-similar diffusive-dispersive approximations for general systems were investigated in LeFloch and Rohde \cite{LeFlochRohde}. Joseph and LeFloch \cite{JLtwo}--\cite{JLfive} extended the technique  
to cover boundary value problems, relaxation approximations, and general diffusion matrices. 
The present paper is the continuation of \cite{JLtwo}--\cite{JLfive}. 
For issues related to the discretization of the viscosity and capillarity terms, we refer 
to \cite{Slemrod-scheme,HayesLeFloch2}.

%======================================================================================== 

\section{Elastodynamics with physical viscosity}
\label{section2}

We begin with the case that $\delta =0$ and $\sigma_w>0$. We will prove: 

\begin{theorem} [Vanishing viscosity limit in elastodynamics]
Assume that the system \eqref{elastic-first} is uniformly hyperbolic, in the sense that 
there exits a positive constant $c_0$ such that 
$$
\inf_{w \in \RR} \sigma_w(w) \geq c_0^2 > 0. 
$$
Then, given arbitrary Riemann data $v_l, w_l$ and $v_r, w_r$, the viscous Riemann problem 
\eqref{elastic2}-\eqref{boundary1} with $\delta =0$ 
admits a solution $v^\eps, w^\eps$ which has uniformly bounded total variation and converges pointwise to a limit $v,w$. 
This limit is a weak solution of the Riemann problem
\eqref{elastic-first}-\eqref{boundary}. 
\end{theorem}

It is worth observing that the solutions $w^\eps, v^\eps$ contain a mild singularity at 
$y=0$: their derivatives up to order $C/\eps$ only are continuous at $y=0$.
Hence, this singularity vanishes in the limit $\eps \to 0$. 
The singularity is due to a factor $y$ appearing in the key equation \eqref{scalar} below, 
in front of the term 
with the highest order derivative. 

\begin{proof}

{\it Step 1. Reduction to a scalar equation for the component $w$. } To simplify the notation we suppress
the subscript $\eps$. 
 We will first study the problem on a finite interval $[-L,L]$ with $L$ chosen to be sufficiently large
and we will show the existence of a solution to \eqref{elastic2} satisfying the boundary conditions
\be
\label{boundary2}
(v,w)(-L) = (v_l, w_l), 
\qquad 
(v,w)(L) = (v_r, w_r).  
\ee

First, we observe that the system can be reduced to a single scalar equation for the unknown $w$.
Namely, using $v' = - y \, w'$ in the first equation of \eqref{elastic2} we find 
\be
\label{scalar}
(y^2 + \eps - \sigma_w(w)) \, w' = - \eps \, y \, w'',
\ee
which will be studied in two regions $[-L,0]$ and $[0,L]$ with the following boundary conditions 
\be
\label{boundl}
w(-L) = w_l, \qquad w(0-) = w_*, 
\ee
\be
\label{boundr}
w(0+) = w_*,  \qquad w(L) = w_r,  
\ee
where $w_*$ is a parameter to be determined so that the boundary condition on the variable $v$
(which has been eliminated from \eqref{scalar}) is satisfied. This condition is obtained by integrating over $[-L,L]$
 the second equation in \eqref{elastic2}: 
\be
\label{boundv}
v_r - v_l = \int_{-L}^L v' \, dy = - L \, w_r + L w_l + \int_{-L}^L w \, dy.  
\ee
This is a global, integral condition on the solution $w$. We will now solve the problem \eqref{scalar}--\eqref{boundv}. 

Our boundary condition can be justified as follows. If we search for a piecewise smooth solution 
of \eqref{elastic2}, then clearly the measure $v' = -y w'$ can not have a mass point at $y=0$, hence 
$v$ is continuous. Similarly, by the first equation in \eqref{elastic2} the function $\sigma(w) - \eps v'$ must be continuous
at $y=0$. But, using again the equation $\eps v' = - y \eps w'$ we see that $\sigma(w)$ itself must be continuous. 
Finally, since $\sigma$ is assumed to be strictly increasing, $w$ must be continuous.

Based on the function $w$ we can define $a(y) : = y - (\sigma_w(w(y)) - \eps)/y$, consider the integral 
$\int_{\alpha^\pm}^y a(x)\, dx$ on $[-L,0)$  and on $(0,L]$ with $\alpha^- \in [-L,0)$  and $\alpha^+ \in (0,L]$, respectively.
Let us show that the minimum with respect to $y$ of this integral is attained at 
$$
\rho^\pm := \pm (\sigma_w(w(\rho^\pm))-\eps )^{1/2}.
$$
For definiteness, consider the case $\alpha=\alpha^- \in [-L,0)$ with $\alpha^- < y$; then
$$
\int_{\alpha}^y a(x)\, dx \geq \int_{\alpha}^y (x - {(\sigma_w^{M}-\eps)\over x})\, dx 
\geq (y^2 - \alpha^2)/2 + \log(\alpha/y)^{\sigma_w^{M} - \eps},
$$
which tends to infinity as $y \rightarrow 0$. Necessarily, the minimum is attained away from $0$. 
As this integral has quadratic growth in $y$,
for $L$ large the minimum $\rho=\rho^-$ is attained in $(-L,0)$ and is given by the equation
$\rho^2 - \sigma_w(w(\rho))+\eps =0$. 

This choice of $\rho^{\pm}$ gives   $\int_{\rho^-}^y a(x) \, dx \geq 0$.
By setting 
$$  
\varphi_-(y) := {e^{-(1/\eps)  \int_{\rho^-}^y a(x) \, dx} \over \int_{-L}^0 e^{-(1 /\eps) \int_{\rho^-}^y a(x) \, dx}}
\qquad 
\varphi_+(y) := {e^{-(1/\eps) \int_{\rho^+}^y a(x) \, dx} \over \int_0^L e^{-(1 /\eps) \int_{\rho^+}^y a(x) \, dx}}, 
$$
the problem under consideration is equivalent to 
\be
\label{equiv}
w(y) = \begin{cases}\displaystyle
w_l + (w_* - w_l) \int_{-L}^y \varphi_- \, dx,     &  y < 0, 
\\\displaystyle
w_r + (w_* - w_r) \int_y^L \varphi_+ \, dx,     &  y > 0, 
\end{cases}
\ee
together with the boundary condition \eqref{boundv}. 

To find $w_*$, we used the continuity of the solution at $y=0$
namely, $w(0-) = w(0+) = w_*$ and  $v(0-) = v(0+) = v_*$.
Integrating the equation $v'= -yw'$ from $-L$ to $0$ and from $0$ to $L$
and using \eqref{equiv}, we get
$$
v_* -v_l = (w_* - w_l)\int_{-L}^{0} -y \varphi_-(y) dy, 
$$
$$
v_r -v_* = (w_* - w_r)\int_{0}^{L} y \varphi_+(y) dy. 
$$
Adding these formulas and solving for $w_*$, we get
\be
\label{wm}
w_* = {v_r - v_l + w_r \int_0^L y \, \varphi_{+} \, dy - 
w_l \int_{-L}^{0} y \, \varphi_{-} \, dy  
\over \int_0^L y \, \varphi_{+}^\delta \, dy -  \int_{-L}^{0} 
y \, 
\varphi_{-}^\delta \, dy}, 
\ee
which provides us with the value of $w_*$. 

Observe that the denominator is bounded away from zero.
To evaluate its minimum value (in terms of the constitutive function $\sigma$ and the Riemann data)
we now need to analyze the functions $\varphi_{\pm}$. 

\

\noindent{\it Step 2. Properties of the functions $\varphi_\pm$.} We will show now that the support of the functions 
$\varphi_\pm$ is essentially concentrated away from the axis $y=0$. 
Let $\sigma_w^m = \min_{-L\leq y \leq L} \sigma_w(w(y))$, $\sigma_w^M = \max_{-L\leq y \leq L} \sigma_w(w(y))$,
$\lambda_\eps^m = (\sigma_w^m-\eps)^{1/2}$ and $\lambda_\eps^M = (\sigma_w^M-\eps)^{1/2}$. 
We claim that there exists a constant
$C>0$ such that 
\be
\label{prop1}
0 < \varphi_{-}(y) \leq {C \over \eps} \begin{cases} 
e^-{(y+\lambda_\eps^M)^2 \over 2\eps}, & -L \leq y < -\lambda_\eps^M,
\\
1, & -\lambda_\eps^M \leq y < -\lambda_\eps^m,
\\ 
e^-{(y+\lambda_\eps^M)^2 \over 2\eps}, & -\lambda_\eps^m \leq y < {-\lambda_\eps^m}/4,
\\
(-2y/\lambda_\eps^m)^{3 \lambda_\eps^m \over 4\eps}, & -{\lambda_\eps^m}/4 \leq y \leq 0,
\end{cases}
\ee
\be
\label{prop2}
0 < \varphi_{+}(y) \leq {C \over \eps} \begin{cases} 
(2y/\lambda_\eps^m)^{3 \lambda_\eps^m \over 4\eps}, & 0 \leq y \leq {\lambda_\eps^m}/4, 
\\
e^-{(y-\lambda_\eps^m)^2 \over 2\eps}, & {\lambda_\eps^m}/4 \leq y < {\lambda_\eps^m},
\\
1, & \lambda_\eps^m \leq y < \lambda_\eps^M,
\\
e^-{(y-\lambda_\eps^M)^2 \over 2\eps}, & \lambda_\eps^M <y <L.
\end{cases}
\ee

This is proved as follows. First, in the interval $0\leq y \leq \lambda_\eps^m/4$ we find 
$$
\aligned
\int_{\rho^+}^y {(x^2 -\sigma_w(x)+\eps)\over x} dx 
&\geq \int_{{\lambda_\eps^m}/2}^y {(x^2 -\sigma_w +\eps)\over x} dx 
\\
&= \int_y^{{\lambda_\eps^m}/2} {(\sigma_w - \eps - x^2)\over x} dx
\\
&\geq  (3/4) (\lambda_\eps^m)^2\int_y^{{\lambda_\eps^m}/2} {dx \over x}
   =\log(\lambda_\eps^m /2 y)^{(3/4) {\lambda_{\eps}^m}^2},
\endaligned
$$
and, next, for $y>\lambda_\eps^M$ 
$$
\aligned
\int_{\rho^+}^y (x^2 -\sigma_w(x)+\eps) dx 
&\geq \int_{\lambda_\eps^M}^y {(x^2 -\sigma_w(x)+\eps) \over x} dx
\\
&= \int_{\lambda_\eps^M}^y (x -(\sigma_w(x)+\eps)^{1/2})(1 +{(\sigma_w(x)+\eps)^{1/2} \over x}) dx
\\
&\geq \int_{\lambda_\eps^M}^y (x-\lambda_\eps^M) dx ={(y - \lambda_\eps^M)^2 \over 2}. 
\endaligned
$$
The other cases are completely similar. Since, on the other hand, it is easy to check that 
$$
 \int_{-L}^0 e^{-(1 /\eps) \int_{\rho^-}^y a(x) \, dx} dy \geq C \eps,
\qquad 
\int_0^L e^{-(1 /\eps) \int_{\rho^+}^y a(x) \, dx} dy \geq C \eps, 
$$
this completes the derivation of the properties of $\varphi_\pm$.

\

\noindent{\it Step 3. Existence result.} 
Observe that $w_*$ consists of the sum of 
$$
{{v_r - v_l}  \over {\int_0^L y \, \varphi_+ \, dy -  \int_{-L}^0 y \, \varphi_- \, dy}}
=  {{v_r - v_l} \over {\int_{-\lambda_\eps^M}^{-\lambda_\eps^m} -y \, \varphi_- \, dy + \int_{\lambda_\eps^m}^{\lambda_\eps^M} y \, \varphi_+ \, dy + O(\eps^n)}}
$$
and a convex combination of $w_l,w_r$.
In view of this observation, together with the estimates \eqref{prop1} 
and \eqref{prop2}, we get  
$$
\aligned 
|w_*| 
& \leq \max (w_l,w_r) + {{|v_r - v_l|} \over {2 \lambda_\eps^m} + O(\eps^n)}
\\
& \leq \max (w_l,w_r) + {{|v_r - v_l|} \over {c_0}} =: \Lambda_0 
\endaligned 
$$
and we arrive at the uniform bound
\be
\label{wstar}
| w_* | \leq \Lambda_0. 
\ee

Replacing now $w_*$ (by its value given in \eqref{wm}) in the formula
\eqref{equiv} we can define a mapping 
$w \in C^0([-L,L]) \mapsto T(w) \in C^0([-L,L])$. We claim that, for fixed $\eps$,  the function $T(w)$ is 
of class $C^1$ and 
$$
\aligned 
& \|T(w)\|_{C^0} \leq \Lambda_0, 
\\
& \|T(w)'\|_{C^0} \leq {C \over \eps} \, (|w_r - w_* | +|w_r - w_* |). 
\endaligned 
$$
Indeed, these estimates follow immmediatly from \eqref{equiv}, \eqref{prop1}, and \eqref{prop2}

In consequence, the operator $T$ is a compact map from the convex bounded set $\big\{ w \, / \,  \|w\|_{C^0([-L,L]} \leq \Lambda_0\big\}$ into
itself. By Schauder's fixed point theorem, $T$ has a fixed point which satisfies 
\eqref{equiv} and, clearly, is of class $C^1$. Furthermore, in view of
\eqref{equiv} and \eqref{wm}, we have 
\be
\label{es1}
\int_{-L}^L |w'(y)| dy \leq | w_r - w_* | +  | w_l - w_* | 
\leq  | w_r - w_l |+ {2 \over c_0} \, | v_r - v_l |
\ee
and
\be
\label{es2}
\int_{-L}^L |v'(y)| dy \leq \int_{-L}^L |y w'(y)| dy 
\leq (\lambda_0^M+1) \, \big( |w_r - w_l|+ {2  \over c_0} \, |v_r - v_l| \big). 
\ee

Next, since all of the estimates are uniform in $L$, we can let $L \to \infty$ and we conclude with 
existence of a solution $(v^\eps,w^\eps)$ on defined on the whole real line $(-\infty,\infty)$. 
This function $(v^\eps,w^\eps)$ is a bounded solution of \eqref{elastic2}-\eqref{boundary1}
and has uniformly bounded total variation
$$
\int_{-\infty}^\infty \big( |w'(y)| + |v'(y)| \big) \, dy \leq (\lambda_0^M+2) \, \big( |w_r - w_l|+ {2 \over c_0} \, |v_r - v_l| \big). 
$$
By Helly's compactness theorem, it admits a subsequence converging pointwise to a limit $(v,w)$, 
as $\eps$ goes to zero. Clearly, 
the limit is a weak solution of the problem \eqref{elastic-first}-\eqref{boundary}.
\end{proof} 

%=============================================================== 

\section{Elastodynamics with physical viscosity and capillarity}
\label{section3}

In this section, we consider the full system with physical viscosity and capillarity included.  We will prove:

\begin{theorem}  [Vanishing viscosity-capillarity limit in elastodynamics]
Assume that $\inf \sigma_w \geq c_0^2 >0$. Then, 
given arbitrary Riemann data $v_l,v_r, w_l, w_r$, the problem 
\eqref{elastic2}-\eqref{boundary1} with $\delta=\gamma \eps^2$,
$\gamma>0$, 
admits a solution $v^\eps, w^\eps$ which has uniformly bounded total variation 
and converges to a limit $v,w$. This limit is a weak solution to the Riemann 
problem \eqref{elastic-first}-\eqref{boundary}. 
\end{theorem}

Observe that this is a large data result. 

\begin{proof}

{\bf Step 1:  reduction to a scalar equation. } 
We will first study the problem on a finite interval $[-L,L]$ with $L$ chosen to be sufficiently large
with boundary conditions 
\be
\label{boundary3}
(v,w)(-L) = (v_l, w_l), 
\qquad 
(v,w)(L) = (v_r, w_r).  
\ee

As before, the system can be reduced to a single scalar equation for the unknown $w$.
Namely, using $v' = - y \, w'$ in the first equation of \eqref{elastic2} we find 
\be
\label{scalar-B}
(y^2 + \eps - \sigma_w(w)) \, w' = - \eps \, y \, w'' - \delta w''',
\ee
which will be studied in two regions $[-L,0]$ and $[0,L]$ with the following boundary conditions 
\be
\label{boundl-BC}
w(-L) = w_l, \qquad w(0-) = w_*, 
\ee
\be
\label{boundr-BC}
w(0+) = w_*,  \qquad w(L) = w_r,  
\ee
where $w_*$ is a parameter to be determined so that the boundary condition on the variable $v$
(which has been eliminated from \eqref{scalar}) is satisfied. This condition is obtained by integrating over $[-L,L]$
 the second equation in \eqref{elastic2}: 
\be
\label{boundv-BC}
v_r - v_l = \int_{-L}^L v' \, dy = - L \, w_r + L w_l + \int_{-L}^L w \, dy.  
\ee
This is a global, integral condition on the solution $w$. We will now solve 
the problem \eqref{scalar-B}--\eqref{boundv-BC}.

We will follow 
a method originally introduced by LeFloch and Rohde \cite{LeFlochRohde} for small amplitude solutions of 
systems of conservation laws.
Setting $\varphi = e^\beta H$, the differential equation under consideration becomes
$$
\delta [H'' + 2 H'\beta' + H \beta'' + H \beta'^2] + \eps y [H' +H \beta'] + (y^2 +\eps - \sigma_w ) H = 0
$$
Setting the coefficient of $H'$ equal to $0$, we get 
$\beta' = {-\eps y \over 2 \delta}$, thus $\beta (y) = {-\eps y^2 \over 4 \delta}$, 
and
\be
H'' = {\mu (w(y),y) \over \delta} H, 
\ee
where
\be
\mu(w(y),y) = \sigma_w + y^2({\eps^2 \over 4 \delta} - 1) - {\eps \over 2}. 
\ee

We take $\delta = \gamma \eps^2$ and, in this case,  
$$
\mu(w,y) =  \sigma_w + y^2({1 \over 4 \gamma} - 1) - {\eps \over 2} > 0, 
$$
provided $\eps,\gamma>0$ are sufficiently small. 
Upon integrating this equation for $H$ and substituting, we get
\be
\phi(y) = { 1+\Phi(y) \over (4\gamma \mu)^{1/4}} e^{p(y,\rho) \over \eps}, 
\ee
where
$$
\aligned 
& |\Phi(y)|+{\eps \gamma^{1/2} \over 2 \mu(y)^{1/2}}|\Phi'(y)| \leq \eps \gamma^{1/2} k, 
\\
& k:=\int_{-L}^L \mu^{-5/4}|\mu'|^2 dx + \int_{-L}^L \mu^{-3/2}|\mu''| dx
\endaligned 
$$
and
\be
\aligned
p(y,\rho) 
&:= \int_\rho^y({-x \over 2\gamma} + \sqrt {\mu(w(x),x) \over \gamma}) dx
\\
&=\int_{\rho}^y ({-x\over 2\gamma} +\sqrt{{\sigma_w -\eps/2 + x^2(1/{4\gamma} -1) \over \gamma}} dx 
\\
&={1\over 2\gamma}\int_{\rho}^y x(-1 + \sqrt{1+4\gamma (\sigma_w -\eps/2 - x^2)} ) dx. 
\endaligned
\ee

Here, $\rho$ is the maximizer of the the integral
$$
 \max_{-L < y < L}\int_\alpha^y({-x \over 2\gamma} + \sqrt {\mu(w(x),x) \over \gamma})dx. 
$$
Indeed $\rho$ satisfies 
$$
\rho^2 - \sigma_w(w(\rho))+\eps/2 =0, 
$$
which has pairs of solutions  $\rho_- \in [-L,0)$  and on $\rho _+ \in (0,L]$,  satisfying 
$$
\rho^\pm = \pm (\sigma_w(w(\rho^\pm))-\eps/2 )^{1/2}.
$$

We consider $\varphi_-$ and $\varphi_+$  on  $[-L,0]$ and $[0, L]$, defined by
\be
\phi_-(y) = { 1+\Phi(y) \over (4\gamma \mu)^{1/4}} e^{p(y,\rho_-) \over 
\eps}
\ee
\be
\phi_+(y) = { 1+\Phi(y) \over (4\gamma \mu)^{1/4}} e^{p(y,\rho_+) \over 
\eps}
\ee
and setting
$$  
\varphi_-(y) := {\phi_-(y) \over \int_{-L}^0 \phi_-(x) \, dx}
\qquad 
\varphi_+(y) := {\phi_+(y) \over \int_0^L \phi_+(x) \, dx}, 
$$
the problem under consideration is equivalent to 
\be
\label{equiv-E}
w(y) = \begin{cases}
w_l + (w_* - w_l) \int_{-L}^y \varphi_- \, dx,     &  y < 0, 
\\
w_r + (w_r - w_*) \int_y^L \varphi_+ \, dx,     &  y > 0, 
\end{cases}
\ee
together with the boundary condition \eqref{boundv}. The latter can be expressed in its original form involving $w'$
(see the second equation in \eqref{elastic2}) and yields 
\be
\label{wm-C}
w_* := {v_r - v_l + w_r \int_0^L y \, \varphi_+ \, dy - w_l \int_{-L}^0 y \, \varphi_- \, dy  
\over \int_0^L y \, \varphi_+ \, dy -  \int_{-L}^0 y \, \varphi_- \, dy}. 
\ee 
Observe that the denominator in the formula for $w_*$ is bounded away from zero.
To find the value of the minimum in terms of the constitutive function $\sigma$ and the Riemann data
we need the following properties of $\varphi_{\pm}$

%__________________________________________

\

\noindent{\bf Step 2: Properties of $\varphi_{\pm}$. } With the same $\sigma_\eps^m$ and $\sigma_\eps^M$
as before,
set $\lambda_\eps^m = (\sigma_w^m-\eps/2)^{1/2}$ and $\lambda_\eps^M = (\sigma_w^M-\eps/2)^{1/2}$.
There exist constants $C_1>0$, and $C>0$ depending only on $\gamma, \sigma_w^m, \sigma_w^M$ 
and such that
\be
\label{prop3}
0 < \varphi_{-}(y) \leq {C_1\over\eps} \begin{cases} 
e^{C(y+\lambda_\eps^M)^{2} \over \eps y}, &  -L \leq y < -\lambda_\eps^M,
\\
1, & -\lambda_\eps^M \leq y < -\lambda_\eps^m,
\\ 
e^-{C(y+\lambda_\eps^m)^2 \over \eps}, & -\lambda_\eps^m \leq y < -\gamma^{1/2}\lambda_\eps^m,
\\
e^-{C|y-\lambda_\eps^m \gamma^{1/2}| \over \eps}, & -\gamma^{1/2}\lambda_\eps^m \leq y \leq 0,
\end{cases}
\ee
\be
\label{prop4}
0 < \varphi_{+}(y) \leq {C_1\over\eps} \begin{cases}  
e^-{C|y-\lambda_\eps^m \gamma^{1/2}| \over \eps}, & 0 \leq y \leq \gamma^{1/2}\lambda_\eps^m,
\\
e^-{C(y-\lambda_\eps^m)^2 \over \eps}, &  \gamma^{1/2} {\lambda_\eps^m} \leq y < {\lambda_\eps^m},
\\
1, & \lambda_\eps^m \leq y < \lambda_\eps^M,
\\
e^-{C(y-\lambda_\eps^M)^2 \over \eps y}, & \lambda_\eps^M <y <L.
\end{cases}
\ee

First, we prove the estimates for $\varphi_+$. Consider the region
$y>\lambda_\eps^M$
$$
\aligned
p(y,\rho_+) 
&=\int_{\rho^+}^y \Big({-x\over 2\gamma} +\sqrt{{\sigma_w -\eps/2 + x^2(1/{4\gamma} -1) \over \gamma}} \Big) \, dx
\\
&\leq \int_{{\lambda_\eps^M}}^y   \Big({-x\over 2\gamma} +\sqrt{{{\sigma_w -\eps/2 + x^2(1/4\gamma -1)} \over \gamma}}\Big) \, dx 
\\
& = {1\over 2\gamma}\int_{{\lambda_\eps^M}}^y x \, \Big(-1 + \sqrt{1+4\gamma {(\sigma_w -\eps/2 - x^2)\over x^2}} \Big) \, dx
\\
&\leq {\lambda_\eps^M\over 2\gamma}\int_{{\lambda_\eps^M}}^y 
      \Big(-1 + \sqrt{1+4\gamma {(\sigma_w -\eps/2 - x^2)\over x^2}} \Big) \, dx,
\endaligned
$$
and, for $y>\lambda_\eps^M$,
$$
\aligned
p(y,\rho_+)
&= \leq {\lambda_\eps^M \over 2\gamma}\int_{{\lambda_\eps^M}}^y
2\gamma {(\sigma_w -\eps/2 - x^2)\over x^2} dx
\\
&\leq {\lambda_\eps^M}\int_{{\lambda_\eps^M}}^y (\sigma_w -\eps/2)^{1/2} - x) {(\sigma_w -\eps/2)^{1/2} +x)\over x^2} dx
\\
&\leq {\lambda_\eps^M}\int_{{\lambda_\eps^M}}^y {(\sigma_w -\eps/2)^{1/2} - x) \over x}dx
\\
&\leq {\lambda_\eps^M \over y}\int_{{\lambda_\eps^M}}^y (\lambda_\eps^M - x)dx
= -{\lambda_\eps^M \over y}{(y -\lambda_\eps^M)^2 \over 2}.  
\endaligned
$$

Now, consider the region  
$\gamma^{1/2}\lambda_\eps^m/2 \leq y \leq \lambda_\eps^m$,
$$
\aligned
p(y,\rho_+) 
&=\int_{\rho^+}^y \Big({-x\over 2\gamma} +\sqrt{{\sigma_w -\eps/2 + x^2(1/{4\gamma} -1) \over \gamma}}\Big) \,  dx
\\
&\leq \int_{{\lambda_\eps^m}}^y \Big({-x\over 2\gamma} +\sqrt{{{\sigma_w -\eps/2 + x^2(1/4\gamma -1)} \over \gamma}}\Big) \, dx 
\\
&\leq \int_y^{\lambda_\eps^m} \Big ({x\over 2\gamma} -\sqrt{{{\sigma_w -\eps/2 + x^2(1/4\gamma -1)} \over \gamma}}\Big) \, dx
\\
&={1\over 2\gamma}\int_y^{{\lambda_\eps^m}} x \, \Big(1 - \sqrt{1+4\gamma {(\sigma_w -\eps/2 - x^2)\over x^2}} \Big) \, dx
\\
&\leq{\lambda_\eps^m\over 4\gamma^{1/2}}\int_y^{{\lambda_\eps^m}}
  \Big(1 - \sqrt{1+4\gamma {(\sigma_w -\eps/2 - x^2)\over \lambda_\eps^m}} \Big) \, dx, 
\endaligned
$$
and, for $\gamma^{1/2}\lambda_\eps^m/2 \leq y \leq \lambda_\eps^m$,
$$
\aligned
p(y,\rho_+) &-\leq {\lambda_\eps^m\over 4\gamma^{1/2}}\int_y^{{\lambda_\eps^m}} 4\gamma {(\sigma_w -\eps/2 - x^2)\over \lambda_\eps^m}) dx
\\
&\leq -C \int_y^{{\lambda_\eps^m}} 4\gamma (\sigma_w -\eps/2)^{1/2} - x) (\sigma_w -\eps/2)^{1/2} +x) dx
\\
&\leq -C \int_y^{{\lambda_\eps^m}} (\sigma_w -\eps/2)^{1/2} - x)dx
\\
&\leq -C \int_y^{{\lambda_\eps^M}} (\lambda_\eps^m - x)dx = -C{(y -\lambda_\eps^m)^2 \over 2}. 
\endaligned
$$

For  $0\leq y \leq \gamma^{1/2}\lambda_\eps^m/2$,
$$
\aligned
p(y,\rho_+) 
&=\int_{\rho^+}^y \Big({-x\over 2\gamma} +\sqrt{{\sigma_w -\eps/2 + x^2(1/{4\gamma} -1) \over \gamma}} \Big) \, dx
\\
&\leq \int_{{\lambda_\eps^m}\gamma^{1/2}}^y
\Big({-x\over 2\gamma} +\sqrt{{{\sigma_w -\eps/2 + x^2(1/4\gamma -1)} \over \gamma}}\Big) \, dx 
\\
&\leq \int_y^{\lambda_\eps^m \gamma^{1/2}} 
\Big({x\over 2\gamma} -\sqrt{{{\sigma_w -\eps/2 + x^2(1/4\gamma -1)} \over \gamma}}\Big) \, dx
\\
&\leq \int_y^{{\lambda_\eps^m}\gamma^{1/2}} ({\gamma^{1/2}\lambda_\eps \over 2 \gamma} - {\lambda_\eps^m \over \gamma^{1/2}})dx
\\
&=\int_y^{{\lambda_\eps^m}\gamma^{1/2}} - {\lambda_\eps^m \over 2\gamma^{1/2}}dx
=-{C|y-\lambda_\eps^m \gamma^{1/2}|}. 
\endaligned
$$
The other cases are identical. One can also check (for some $c>0$) 
$$
 \int_{-L}^0 \phi_-(y) dy \geq c \, \eps,
\qquad 
\int_0^L \phi_+(y) dy \geq c \, \eps, 
$$
so that the desired properties of $\varphi_+,\varphi_-$, follow.

\

\noindent{\bf Existence arguments. }
Now $w_*$ is sum of convex combination of $w_l,w_r$ and the quantily
$$
 {{v_r - v_l}  \over {\int_0^L y \, \varphi_+ \, dy -  \int_{-L}^0 y \, \varphi_- \, dy}}
=  {{v_r - v_l} \over {\int_{-\lambda_\eps^M}^{-\lambda_\eps^m} -y \,
\varphi_- \, dy + \int_{\lambda_\eps^m}^{\lambda_\eps^M} y \, \varphi_+ \, dy + O(\eps^n)}}.
$$
From this we get 
% for any $n$,
$$
|w_*| \leq \max (w_l,w_r) + {{|v_r - v_l|} \over {2 \lambda_\eps^m} + O(\eps^n)}\leq \Lambda_0. 
$$
We can now replace $w_*$ (given by \eqref{wm-C}) in \eqref{equiv-E} and
arrive at a mapping 
$w \in C^0([-L,L]) \mapsto T(w) \in C^0([-L,L])$. For fixed $\eps$,  $T(w)$ is of class $C^1$, and 
$$
\aligned 
& \|T(w)\|_{C^0} \leq \Lambda_0, 
\\
& \|T(w)'\|_{C^0} \leq {C \over \eps} \, (|w_r - w_* | +|w_r - w_* |). 
\endaligned 
$$
These estimates follow from \eqref{equiv-E}  \eqref{prop3} and
\eqref{prop4}

Thus, $T$ is a compact map from a convex bounded set $\{ w : \|w\|_{C^0} \leq \Lambda_0\}$ into
itself. By Schauder's fixed point theorem, $T$ has a fixed point in $C^0[-L,L]$ and satisfies 
\eqref{equiv-E} and so the solution is of class $C^1$. Furthermore, in
view of \eqref{equiv-E} and \eqref{wm-C}, 
\be
\label{es1-B}
\int_{-L}^L |w'(y)| dy \leq | w_r - w_* | +  | w_l - w_* | \leq  | w_r - w_l |+ 2{| v_r - v_l | \over c_0}
\ee
and
\be
\label{es2-B}
\int_{-L}^L |v'(y)| dy \leq \int_{-L}^L |y w'(y)| dy \leq (\lambda_0^M+1) \, \Big( |w_r - w_l|+ 2{|v_r - v_l| \over c_0}\Big). 
\ee

As the estimates are uniform in $L$, we have existence of solution $(v^\eps,w^\eps)$ on $(-\infty,\infty)$
we have $(v^\eps,w^\eps)$ bounded solution of \eqref{elastic2} and \eqref{boundary1} with uniform total variation
$$
\int_{-\infty}^\infty (|w'(y)|+|v'(y)| dy \leq (\lambda_0^M+2)[ |w_r - w_l|+ 2{|v_r - v_l| \over c_0}]
$$
By compactness there exists a sequence converges in $L^1$ and pointwise
a.e. to a function $(v,w)$ as $\eps$ goes to zero and solves
\eqref{elastic-first} and \eqref{boundary} with $\eps =0$. 
\end{proof}

%=======================================================================================

\section{Phase dynamics with physical viscosity}
\label{section4}

We turn to the model of phase transition dynamics when the system is not strictly
hyperbolic, that is, when $\sigma_w(w)$ is only non-negative, or even is hyperbolic-elliptic 
when $\sigma_w$ takes negative values. 

To illustrate the key difficulty we will have to cope with, let us first consider the example
$$
a(y) = y + c/y, 
$$
with $c$ a constant. When $c\geq 0$, $\rho_{+} =\delta$ and 
$\int_{\rho^+}^y a(x) dx = {y^2-\rho_{+}^2 \over 2} + 
\log(|y/\rho_{+}|^{c})$ and so 
$$
e^{-{1\over\eps}\int_{\rho_{+}}^y a(x) dx} = e^{-(y^2-\rho_{+}^2) \over 
2 \eps} (y/\rho_{+})^{-c/\eps} = e^{-(y^2-\delta^2) \over 
2 \eps} (y/\delta)^{-c/\eps}.
$$
In this case 
\[
\varphi_+(y) =  {e^{-(y^2-\delta^2)\over 
2 \eps} (y/\delta)^{-c/\eps}\over \int_{\rho_{+}}^L e^{-(y^2-\delta^2)
\over 2 \eps} (y/\delta)^{-c/\eps} dy}, 
\]
which is concentrated at $y=\delta$. For $c<0$, we have, 
$\rho_+ = (- c)^{1/2}$ and
$$
e^{-{1\over\eps}\int_{\rho_{+}}^y a(x) dx} = e^{-(y^2+c) \over 
2 \eps} (y/(-c)^{1/2})^{-c/\eps},
$$
so that
\[
\varphi_+(y) =  {e^{-(y^2+c) \over 
2 \eps} (y/(-c)^{1/2})^{-c/\eps} \over \int_{\rho_{+}}^L
e^{-(y^2+c)\over 2 \eps} (y/(-c)^{1/2})^{-c/\eps} dy}, 
\]
which is concentrated at $y=(-c)^{1/2}$.

In the previous section, to show that $w_*$ is bounded independent of
$\eps$ it was crucial that $\varphi^{+}$ and $\varphi^{-}$ be concentrated
away from $0$. In the above example, if $c \geq 0$, we have seen
that this is not the case. This is the difficulty is getting estimates when $\sigma_w(w(y))$
oscillates, i.e.~takes both negative and positive values.

We now prove the following result: 

\begin{theorem}  [Vanishing viscosity limit in phase dynamics]
Suppose that the first order system \eqref{elastic-first} is uniformly
hyperbolic for large value of $w$ but
may admit elliptic regions in the phase space, that is for some constants $M,c_0 >0$
$$
\inf_{|w| \geq M} \sigma_w(w) \geq {c_0}^2. 
$$
Suppose also that there exist constants $c_1>0$ and $\eta>0$ such that
$|\sigma_w (w) |\leq c_1 |w|^{2-\eta}$ for $|w|>1$.  
Then for all Riemann data $v_l, w_l$ and $v_r,w_r$ in the hyperbolic region of the phase space, i.e. 
$$
w_l, w_r \in \Hcal := \big\{w  \, / \, \sigma_w >0 \big\}
$$
the viscous Riemann problem \eqref{elastic2}-\eqref{boundary1},  
with $\delta = 0$, 
admits a solution $v^\eps, w^\eps$ which has uniformly bounded total variation 
at least away from $y=0$, and more precisely
$$
\int_\RR \big( |v'_\eps| + y \, |w_\eps'| \big) \, dy \lesssim |v_r - v_l| + | w_r - w_l|. 
$$
The functions $v_\eps, w_\eps$ converge pointwise at all $y \neq 0$ to a limit $v,w$, 
which is a solution of the Riemann problem
\eqref{elastic-first}-\eqref{boundary} {\rm away from the axis} 
$y=x/t =0$. Furthermore, $v$ has bounded variation and so 
admits left- and right-hand limits at $y=0$, while the function $w$ and its variation measure $dw/dy$ satisfy 
$$
\aligned 
& |w| \lesssim {1 \over |y|}, \quad y \neq 0,  
\\
& \int_\RR |y| \, |{dw \over dy}| < \infty. 
\endaligned 
$$
\end{theorem} 

Hence, the component $v$ only has globally bounded variation. 
In turn, the conservation laws \eqref{elastic-first} are satisfied in the two regions $x < 0$
and $x>0$, but a singularity may arise on the axis.

\begin{proof}
{\it Step 1. A priori estimates for the components $v,w$.} To simplify the notation we suppress the subscript $\eps$. 
We will first study the problem away from the axis $y=0$. We consider the function 
\be
a(y) = y - (\sigma_w(w(y)) - \eps)/y,   
\ee
and, given some $\delta>0$, we study the problem  in two regions $[-L,-\delta]$ and 
$[\delta,L]$, with the following boundary conditions 
\be
\label{boundl-B}
w(-L) = w_l, \qquad w(\delta-) = w_*, 
\ee
\be
\label{boundr-B}
w(\delta) = w_*,  \qquad w(L) = w_r.  
\ee
We set 
$$  
\varphi_{-}^\delta(y) := {e^{-(1/\eps) \int_{\rho_-}^y a^\delta(x) \, 
dx} \over \int_{-L}^{-\delta} e^{-(1 /\eps) \int_{\rho_-}^y a^\delta(x) 
\, dx}}, 
\qquad 
\varphi_{+}^\delta(y) := {e^{-(1/\eps) \int_{\rho_+}^y a^\delta(x) \, 
dx} \over \int_\delta^L e^{-(1 /\eps) \int_{\rho_+}^y a^\delta(x) \, 
dx}}, 
$$
on $[-L,-\delta]$ and $[\delta,L]$
where $\rho_\pm$ are the points where $\int^y a(x) dx$ attains its global minimum.
 
Let us consider  $\lambda_\eps^{M+} = \sup_{1 \leq y \leq 
L}{(\sigma_w(w)-\eps)^+}^{1/2}$ 
and let $c> \max\{1,\lambda_\eps^{M+}\}$. For $y>c$
$$
\aligned
\int_{\rho^+}^y (x^2 -\sigma_w(x)+\eps) dx 
&\geq \int_{c}^y {(x^2 -\sigma_w(x)+\eps) \over x} dx
\\
&= \int_{c}^y ((x -(\sigma_w(x)+\eps)^{+})^{1/2})(1 
+{((\sigma_w(x)+\eps)^{+})^{1/2} \over x}) dx
\\
&\geq \int_{c}^y (x-c) dx = {(y - c)^2 \over 2}. 
\endaligned
$$
So, for $y>(1+ \sup_{1 \leq y \leq L}({(\sigma_w(w)-\eps)^+})^{1/2}$, we get 
\be
\label{prop6}
\varphi(y) dy \leq {1\over\eps}e^{-{(y - c)^2 \over 2\eps}}
\ee
and
\be
\label{wavesupport}
\int_\delta^L y\varphi(y) dy \leq 
(2+ \sup_{1 \leq y \leq 
L}{((\sigma_w(w)-\eps)^+})^{1/2}. 
\ee 
We then set
\be
\label{equiv2}
w(y) = \begin{cases}
w_l + (w_* - w_l) \int_{-L}^y \varphi_{-}^\delta \, dx,     &  y < 
-\delta, 
\\
w_r + (w_* - w_r) \int_y^L \varphi_{+}^\delta \, dx,     &  y >\delta.
\end{cases}
\ee

Here, we have taken $w(-\delta) = w(\delta) = w_*$ and we also take 
$v(-\delta) = v(\delta) = v_*$. Integrating the equation $v'=yw'$ from $-L$ to $-\delta$ and 
$\delta$ to $L$  and using \eqref{equiv2}, we get
$$
v_* -v_l = (w_* - w_l)\int_{-L}^{-\delta} -y \phi_-(y) dy
$$
and 
$$
v_r -v_* = (w_* - w_r)\int_{\delta}^{L} y \phi_+(y) dy. 
$$
Adding these formulas, we get
\be
\label{wm2}
w_* := {v_r - v_l + w_r \int_\delta^L y \, \varphi_{+}^\delta \, dy - 
w_l \int_{-L}^{-\delta} y \, \varphi_{-}^\delta \, dy  
\over \int_\delta^L y \, \varphi_{+}^\delta \, dy -  \int_{-L}^{-\delta} 
y \, 
\varphi_{-}^\delta \, dy} 
\ee
and, subtracting,
\be
\label{vm2}
v_* = v_l +{[v_r - v_l + (w_r - w_l) \int_\delta^L y \, 
\varphi_{+}^\delta \, dy]
 \int_{-L}^{-\delta} -y \, \varphi_{-}^\delta \, dy  
\over \int_\delta^L y \, \varphi_{+}^\delta \, dy -  \int_{-L}^{-\delta} 
y \, \varphi_{-}^\delta \, dy}. 
\ee

Now, $v$ can be expressed in terms of the boundary data $w_l,w_r,v_l,v_r$ and $\delta$: 
\be
\label{equiv3}
v^\delta(y) = \begin{cases}
 v_l +{[v_r - v_l + (w_r - w_l) \int_\delta^L y \, \varphi_{+}^\delta \, 
dy] \int_{-L}^y -x \, \varphi_{-}^\delta \, dx  
\over \int_\delta^L y \, \varphi_{+}^\delta \, dy -  \int_{-L}^{-\delta} 
-y \,  \varphi_{-}^\delta \, dy}, & -L \leq y<-\delta
\\
 v_r + {[v_l - v_r + (w_r - w_l) \int_-L^\delta -y \, \varphi_{-}^\delta 
\, dy] \int_{y}^L x \, \varphi_{+}^\delta \, dx  
\over \int_\delta^L y \, \varphi_{+}^\delta \, dy +  \int_{-L}^{-\delta} 
\, -y\varphi_{-}^\delta \, dy}.     & L> y > \delta. 
\end{cases}
\ee
Also, $w(y)$ can be written as
\be
\label{equiv4}
w^\delta(y) = \begin{cases} \displaystyle 
w_l + (v_* - v_l){ \int_{-L}^y \varphi_{-}^\delta \, dx, \over 
\int_{-L}^{-\delta} y\varphi_{-}^\delta \, dx}  &  y < -\delta, 
\\  \displaystyle
w_r + (v_r - v_*){\int_y^L \varphi_{+}^\delta \, dx, \over \int_\delta^L 
y\varphi_{+}^\delta \, dx}  &  y > \delta.
\end{cases}
\ee

Now clearly $v_{*}$ is bounded independent of $\eps>0$,$\delta>0$ :
$$
|v_* -v_l|\leq |v_r - v_l| + L (|w_r -w_l|) 
$$
and from \eqref{equiv3} and \eqref{equiv4}, we get the following estimates
$$
\aligned 
& |y| \, |w^\delta(y)| \leq (|v_r| + |v_l|) + 2L(|w_r| +|w_l|),
\\
& |v^\delta(y)|\leq 2 (|v_r| + |v_l|) + L (|w_r| +|w_l|), 
\endaligned 
$$
Also we have the following estimates on the derivatives of
$v^\delta,w^\delta$.
$$
\aligned 
& \int_{L>|y|>\delta}|y w^\delta(y)'| dy \leq {(|v_r - v_l|) +
L(|w_r-w_l|)},
\\
& \int_{L>|y|>\delta}|v^\delta(y)'| dy\leq (|v_r-v_l|) + L (|w_r-w_l|). 
\endaligned 
$$ 
Furthermore, for all $L>|y|>\delta$
$$
\aligned 
& |v^\delta(y)'| \leq {1  \over \eps \delta} \Big( (|v_r-v_l|) + L (|w_r-w_l|)y \Big), 
\\
& |w^\delta(y)'| \leq {1 \over \eps \delta} \Big( (|v_r-v_l|) + L (|w_r-w_l|)\Big). 
\endaligned 
$$

The existence of solution for fixed $L$ follows from these estimates
as we will see. In order to pass $L\rightarrow \infty$, we need estimates
independent of $L$. For this
we use the growth condition on $\sigma_w$. From the expression
\eqref{vm2} for $v_{*}$, and \eqref{wavesupport}, we get

$$
|v_* -v_l|\leq |v_r - v_l| + (2+ \sup_{1 \leq y \leq 
L}{(\sigma_w(w)-\eps)^+}^{1/2}) |w_r -w_l|. 
$$

Using this in \eqref{equiv3} and \eqref{equiv4}, we get
$$
\aligned 
&\sup_{1\leq y\leq L} |w^\delta(y)| \leq (|v_r| + |v_l|) + ((3+ \sup_{1
\leq y \leq L}{(\sigma_w(w)-\eps)^+}^{1/2})(|w_r| +|w_l|),
\\
& |v^\delta(y)|\leq 2 (|v_r| + |v_l|) + ((2+ \sup_{1
\leq y \leq L}{(\sigma_w(w)-\eps)^+}^{1/2})(|w_r| +|w_l|), 
\endaligned 
$$   
Now using our assumption
$|\sigma_w|\leq c_1 |w|^{2-\eta}$ for $|w|>1$ in the first inequality, we
conclude that $\sup_{1\leq y\leq L}|w^\delta(y)|$ is independent of $L$. Similar
estimate holds for $[-L\leq y\leq-1]$. We get there exists a constant
$C=C(v_l,v_r,w_l,w_r)$, independent of $L$ such that
\be
\label{bound[L]}
\sup_{1\leq |y|\leq L}|w^\delta(y)| \leq C
\ee
With this constant $C$, let $\lambda^M =sup_{|w|\leq
C}{|\sigma_w(w)|}^{1/2}$.
From \eqref{prop6},
we get  $\varphi^{\delta}_{+}$ is essentially supported in 
$[\delta, \lambda^M+1]$. 
Similar arguments give $\varphi^{\delta}_{-}$ is essentially
supported in the interval $[-\lambda^M -1,-\delta]$.
It follows from \eqref{equiv3} and \eqref{equiv4} that
there exists a constant $C_1$ 
\be
\label{var}
\aligned 
& \int_{|y|>\delta}|y w^\delta(y)'| dy \leq {(|v_r - v_l|) +
C_1(|w_r-w_l|)},
\\
& \int_{|y|>\delta}|v^\delta(y)'| dy\leq (|v_r-v_l|) + C_1 (|w_r-w_l|). 
\endaligned 
\ee

\noindent{\it Step 2. Existence proof. }
For each $\delta>0$, we can apply Schauder's fixed point theorem, as was explained earlier, 
and we obtain a solution $(w^\delta,v^\delta)$ defined on the interval $[-L,-\delta] \bigcup [\delta, L]$
and satisfying uniform in $\eps>0$ total variation estimate. 
The estimates \eqref{bound[L]} and \eqref{var} allows as to let $L$ tend
to infinity.
By compactness, we have a solution $(v,w)$ 
defined in the region $|y|>\delta$. Since $\delta>0$ is arbitrary, we obtain a well-defined 
solution away from $y=0$. Indeed, $v^\delta$ has a uniform total variation on $R$ and 
$v$ admits left- and right-limit at  $y=0$.
One can check easily that the limit is a weak solution of the problem \eqref{elastic-first}-\eqref{boundary}
away from the axis $y=0$ at least.
Without a control on $\sigma(w)$ near $y=0$ we can not exclude a
concentration term on the axis. 
\end{proof}

%====================================================================================

\section{Boundary Riemann problem and further generalizations}
\label{section5}

In this section, we outline how the method developed in previous sections can be
generalized to the boundary-value problem and indicate several generalizations.
It is well-known that in the strictly hyperbolic case with nondegenerate function $\sigma$, 
the boundary Riemann problem (in $x>0, t>0$) 
\be
\label{elastic51}
v_t -\sigma(w)_x =0, 
\qquad 
w_t - v_x =0.  
\ee
with initial and boundary conditions
\be
\label{boundary51}
w(0+,t) = w_b, 
\qquad 
v(x,0) = v_r, w(x,0) = w_r
\ee
is well posed. This easily follows from an analysis of the wave curves for the system corresponding to the left-moving 
and right-moving characteristic families. The physical regularizations
considered in earlier sections is well-suited to handle this boundary value problem, without producing
boundary layers in the limit. 

Consider the nonlinear elastodynamics with physical viscosity
\be
\label{elastic52}
-y v'- \sigma_w(w)) \, w' = \eps \, v'',
\qquad
-y \, w'-v'=0, \qquad \text{ on $[0,\infty)$, }
\ee
with boundary conditions 
\be
\label{boundary52}
w(0) = w_b, 
\qquad 
v(\infty) = v_r, w(\infty) = w_r.  
\ee
Observe that, for the "full'' viscosity approximation, we need to
prescribe the component $v$ at $y=0$ and this generates a boundary layer
at $y=0$, after passage to the limit $\epsilon \to 0$; see~\cite{JLthree,JLfour}. In
the present case, we show that no boundary layer arises.

As in Section 2,  the problem can be reduced to a scalar equation
for the unknown $w$ 
\be
\label{scalar51}
(y^2 + \eps - \sigma_w(w)) \, w' = - \eps \, y \, w'', \qquad \text{ on $[0,L]$,} 
\ee
for sufficiently large $L$, with boundary conditions,
\be
\label{boundary53}
w(0) = w_b, w(L) = w_r.  
\ee
Once we have $w$, we get  the component $v$ from the equation $yw'+v' =0$ and the boundary condition $v(L)=v_r$. 

The fixed point argument in Section 2 yields a solution $w^\eps$ of 
\eqref{scalar51}-\eqref{boundary53}
which is of uniformly bounded variation. Furthermore, it can be
represented by an integral formula in terms of the function $\varphi_+$
introduced in Section 2 
\[
w^\eps(y) = w_r + (w_b - w_r) \int_y^L \varphi_+ \, dx,    \qquad y > 0. 
\]
Using that $\int_0^L \varphi_+(x) dx =1$, we get
\be
\label{equiv51}
w^\eps(y) = w_b + (w_b - w_r) \int_0^y \varphi_+ \, dx,    \qquad y > 0. 
\ee
Using \eqref{equiv51} in the equation $v'=-yw'$ and integrating from $y$
to
$L$, we get
\be
\label{equiv52}
v^\eps(y) = v_r + (w_b - w_r) \int_y^L x \varphi_+ \, dx,   \qquad   y > 0, 
\ee
where we used $v(L)=v_r$. Note that 
we cannot prescribe boundary condition at $v(0)$ arbitrarily, since from the above equation it follows that 
\[
v(0) = v_r + (w_b - w_r) \int_0^L x \varphi_+ \, dx.
\]

As in Section 2, using the fact that $\varphi_+$ decays exponentially
we can let $L$ tend to $\infty$. Again,
using the properties \eqref{prop2} of $\varphi_+$ near $y=0$ in  
\eqref{equiv51}, it follows easily that there exists $C>0$, a constant independent of $\epsilon$, such
that for all $y>0$ close to the origin 
\[
|w(y) - w_b| \leq C \, y, 
\]
Hence, no boundary layer arises in this approximation at $y=0$ and in the
limit as $\epsilon \rightarrow 0$, the limit function $w=w(y)$
satisfies the boundary condition $w(0+)=w_b$. We summarize our results in
the following theorem:

\begin{theorem}
Assume that $\inf \sigma_w \geq c_0^2 >0$. Then
given arbitrary Riemann boundary data, the problem
\eqref{elastic52}-\eqref{boundary52} 
admits a solution $v^\eps, w^\eps$ which has uniformly bounded total variation 
and converges to a limit $v,w$ which is a solution of the boundary Riemann
problem 
\eqref{elastic51}-\eqref{boundary51}. 
\end{theorem}

\ 

We now discuss a general system. The basic nature of the
system \eqref{elastic51} is that $v$ appear linearly and the
characteristic speeds are equal in magnitude and opposite in sign.
Further using the second equation, the system can be reduced to uncoupled
equations for $w$. After this reduction we
solved for $w$ first and then for $v$. Our results can be
generalized for systems of first order equations for $w$ and $v$ vector
valued functions, having same structure, using the ideas of the
present work and the earlier work \cite{JLthree,JLfour} on boundary value problems. 

For example
let $F : R^N \rightarrow R^N$ be a smooth function with $A(w)=D_{w}F$
has real distinct positive eigenvalues $\lambda_1(w)^2 < \lambda_2(w)^2
<...< \lambda_N(w)^2$
with a complete set of right-left normalized eigenvectors $r_j,
j=1,2,...N$, $l_j,j=1,2,...,N$, $l_j.r_k =\delta_{jk}$. 

We assume that  there exists $c_0>0$ such that $\lambda_j(w) \geq c_0$ for all 
$j=1, 2, ...N$ and for all $w \in B(\delta_0)$, a ball of radius $\delta_0$ around a fixed state
that we take to be $0$. We consider the system of $2N$ equations
\be
\label{general51}
\aligned 
& {\del v \over \del t} - {\del \over \del x} \Big( F(w)) = 0,   
\\
& {\del w \over \del t} - {\del v \over \del x} = 0,
\endaligned 
\ee
where $v=v(t,x)$ and $w=w(t,x)$ are $R^N$ valued functions. The system
\eqref{elastic51} is a special case of \eqref{general51} with $N=1$ and
$\lambda(w)^2=\sigma_w(w)$

We consider boundary value problem for \eqref{general51} on $x>0, t>0$.
When physical viscosity terms are taken into account, 
the corresponding set of differential equations becomes
\be
\label{general52}
\aligned 
& -y \, v' - F(w)' = \eps \, v'', 
\\
& - y \, w' - v' = 0,
\endaligned 
\ee
supplemented with the boundary conditions 
\be
\label{boundary54}
\lim_{y \to \infty} (v,w)(y) = (v_r, w_r), 
\qquad 
w(0+) = w_b.  
\ee

We first consider \eqref{general52} on $[0,L]$ with boundary conditions
\be
\label{boundary55}
(v,w)(L) = (v_r, w_r), 
\qquad 
w(0+) = w_b.  
\ee
As in the previous case the problem can be reduced to first solving a
system for $w(y)$ namely
\be
\label{scalar33}
(y^2 + \eps - A(w)) \, w' = - \eps \, y \, w'',
\ee
on $[0,L]$, with boundary conditions,
\be
\label{boundary7}
w(0) = w_b, w(L) = w_r.  
\ee
Existence of uniformly BV solutions of this problem easily
follows from the work of Joseph and LeFloch \cite{JLthree,JLfour}. We just
outline
the main steps omitting the details. 

We decompose $u'$ in the basis of eigenvectors of $A(u)$, 
\be
\label{expan51}
u'(y)=\sum_{j=1}^{N}a_j r_j(u), a_j= \la l_j,  u' \ra   
\ee
A straight forward calculation lead to the following system of nonliear
equations for $a_j, j=1,2, ..., N$
\be
\label{coupled51}
\eps y a_j' + (y^2 +\eps -\lambda_j^2(w))a_j = D_1(a,a)
\ee
where
\be
\label{quadra51}
D_1(a,a) = -\eps y \sum_{k,i} a_i a_k (D_u r_i.r_k)
\ee

Observe that the linearized form of the equation \eqref{coupled51} has the
form
\[
\eps \, y \, a_j'(y^2 + \eps - \lambda_j^2(w)) \, a_j' = 0,
\]
and the corresponding wave measure $\varphi_{j+}$ has exactly same
properties
\eqref{prop2} given in section 2. 
These wave measures naturally appear when we invert the linear part of
the equation \eqref{coupled51}. Because of the quadratic righthand side
this leads to interaction terms with different families of wave measures
which is controlled by $\sum_{j=1}^N \varphi_{j+}$. Fixed point arguments 
give a uniform BV
solution $w^\eps(y)$ for the system \eqref{scalar33} and
\eqref{boundary7}.
The details are similar to Joseph and LeFloch \cite{JLthree} and, therefore, are omitted. Then 
the BV estimate for $v$ follows from the second equation in 
\eqref{general52} and
\[
v(y)=v_r + \int_y^L x w'(x) dx.
\]

As the estmate \eqref{prop2} shows that
$\varphi_{j+}$ is essentially supported in $[\lambda_j^m,\lambda_j^M]$
we can let $L$ go to infinity. Here $\lambda_j^m$ and $\lambda_j^M$
denotes the minimum and maximum of $\lambda_j(w)$ on the ball
$B(\delta_0)$. We have the following result. 

\begin{theorem} 
Assume that $\inf \lambda_j(w) \geq c_0 >0$. Then there exists
$\delta_1>0,\delta_2>0$ such that, for $w_b, w_r \in B(\delta_1)$
and arbitrary $v_r$ 
the boundary Riemann problem \eqref{general52}-\eqref{boundary54} 
admits a solution $v^\eps, w^\eps$ which has uniformly bounded total variation
with $w^\eps(y) \in B(\delta_2)$  
and converges to a limit $v,w$ which is a solution of the equation
\eqref{general51} with boundary conditions
$w(0+)=w_b, w(x,0)=w_r, v(x,0)=v_r$. 
\end{theorem}

\begin{remark} This analysis is easily extended to the case with physical
viscosity and capillarity as well namely
\be
\label{general53}
\aligned 
& -y \, v' - F(w)' = \eps \, v'' - \gamma \eps^2 w''', \gamma >0, 
\\
& - y \, w' - v' = 0,
\endaligned 
\ee
supplemented with the boundary conditions 
\be
\label{boundary56}
\lim_{y \to \infty} (v,w)(y) = (v_r, w_r), 
\qquad 
w(0+) = w_b.  
\ee
\end{remark}

As in the previous case, the problem can be reduced to first solving a
system for $w(y)$, namely
\be
\label{scalar53}
(y^2 + \eps - A(w)) \, w' = - \eps \, y \, w''-\gamma \eps^2 w''',
\ee
on $[0,L]$, for sufficiently large $L$, with boundary conditions 
\be
\label{boundary57}w(0) = w_b, w(L) = w_r.  
\ee
We decompose $u'$ in the basis of eigenvectors of $A(u)$, 
\be
\label{expan}
u'(y)=\sum_{j=1}^{N}a_j r_j(u), a_j = \la l_j, u'\ra   
\ee
 and then the nonlinear
equations for $a_j, j=1,2, ..., N$ take the form
\be
\label{coupled52}
\gamma \eps^2 a_j'' +\eps y a_j' + (y^2 +\eps -\lambda_j^2(w))a_j =
D_1(a,a)
+ D_2(a,a) +D_3(a,a,a), 
\ee
where
\[
D_1(a,a) = -\eps y \sum_{k,i} a_i a_k (D_u r_i.r_k)
\]
\[
D_2(a,a) = -\gamma \eps^2 (\sum_{k,i} a_i a_k' (D_u r_i.r_k)+
          \sum_{k,i} (a_i a_k)' (D_u r_i.r_k
\]
\[
D_3(a,a,a) = -\gamma \eps^2  \sum_{k,i,l} a_i a_k a_l D_u(D_u r_i.r_k)r_l
\] 

The linearized equation for \eqref{coupled52} is
\[
\gamma \eps^2 a_j'' +\eps y a_j' + (y^2 +\eps -\lambda_j^2(w))a_j = 0, 
\]
and hence the wave measures $\varphi_{j+}$ in this case has the same
qualitative
properties as \eqref{prop4}.

Existence of uniformly BV solutions of the problem \eqref{scalar53} and
\eqref{boundary57} can be easily deduced from \cite{LeFlochRohde}. Then, the BV estimate for $v$ follows from
the second equation of \eqref{general53}. The details are omitted. We get the following result: 

\begin{theorem} 
Assume that $\inf \lambda_j(w) \geq c_0 >0$. Then there exists
$\delta_1$
and $\delta_2$ such that for $w_b, w_r \in B(\delta_1)$ and for arbitrary
$v_r$,  
the problem \eqref{general53}-\eqref{boundary56} 
admits a solution $v^\eps, w^\eps$, which has
uniformly bounded total variation 
with $w^\eps(y) \in B(\delta_2)$
and converges to a limit $v,w$ which is a solution of the equation
 \eqref{general51} with boundary conditions
$w(0+)=w_b,w(x,0)=w_r,v(x,0)=v_r$. 
\end{theorem}

Now let us consider the phase transition case with physical viscosity
and capillarity. Following
previous sections , the system can be reduced to a
single scalar equation for the unknown $w$, 
\be
\label{scalar-B1}
(y^2 + \eps - \sigma_w(w)) \, w' = - \eps \, y \, w'' - \gamma \eps^2 
w''', \gamma>0
\ee
which can be studied in two regions $[-L,-\delta_0]$ and $[-\delta_0,L]$ 
away from $y=0$, with the boundary conditions 
\be
\label{boundl-BC1}
w(-L) = w_l, \qquad w(-\delta_0) = w_*, 
\ee
\be
\label{boundr-BC1}
w(\delta_0) = w_*,  \qquad w(L) = w_r,  
\ee
where $w_*$ is given by \eqref{wm2} with $\delta$ replaced by $\delta_0$ chosen as follows. 
Following Sections 3 and 4, we can construct a BV solution $w$ for 
\eqref{scalar-B1}, in $[-L,-\delta_0]$ and $[\delta_0, L]$
with the boundary conditions \eqref{boundl-BC1} \eqref{boundr-BC1}, provided  
$$
\mu(w,y) =  \sigma_w + y^2({1 \over 4 \gamma} - 1) - {\eps \over 2} > 0. 
$$
This is the case if we impose the restriction $|y|>\delta_0$ where 
\be
\label{delta}
\delta_0> ({4 c \gamma \over 1-4\gamma})^{1/2}. 
\ee
Here, $c>0$ is a constant such that $\sigma_w(w)\geq -c$ for all $w$.

We get the following result : 

\begin{theorem}  [Vanishing viscosity-capillarity limit in phase dynamics]
Suppose that the first-order system \eqref{elastic-first} is uniformly
hyperbolic for large value of $w$ but may admit elliptic regions in the phase space, 
that is for some constants $M,c_0 >0$
$$
\inf_{|w| \geq M} \sigma_w(w) \geq {c_0}^2. 
$$
Suppose also that there exist positive constants $c_1, \eta$ such that
$|\sigma_w|\leq c_1 |w|^{2-\eta}$ for $|w|>1$. 
Let $c>0$ be such that $\sigma_w(w)>-c$ for all $w$ and suppose $\delta_0$ is chosen as in \eqref{delta}.
Then, for all Riemann data $v_l, w_l$ and $v_r,w_r$ in the hyperbolic region of the phase space, i.e. 
$$
w_l, w_r \in \Hcal := \big\{w  \, / \, \sigma_w >0 \big\}
$$
the viscous-capillarity Riemann problem
\eqref{elastic2}-\eqref{boundary1}, with $\delta=\gamma \eps^2$, $\gamma>0$   
admits a solution $v^\eps, w^\eps$ which has uniformly bounded total variation 
in $|y|>\delta_0$, and more precisely
$$
\int_{|y|>\delta_0} \big( |v'_\eps| +  |w_\eps'| \big) \, dy \lesssim
|v_r - v_l| + | w_r - w_l|. 
$$
The functions $v_\eps, w_\eps$ converge pointwise at all $|y|>\delta_0$ to
a limit $v,w$, which is a solution of the Riemann problem
\eqref{elastic-first}-\eqref{boundary} for  $|y|=|x/t| >\delta_0$.
\end{theorem} 

For instance, the example $\sigma(w)=u(u^2-1)$ satisfies all the assumptions of the theorem.
We can not exclude that the vanishing viscosity-capillarity approximations to the Riemann problem 
could contain a highly oscillating (i.e.~weakly but not strongly converging) pattern 
near the axis; its thickness would be $2\delta_0$, at most. The existence of such stationary waves 
is not surprising as it has been observed numerically; they were never pointed out analytically until now. 

%*********************************************************************************************************

\section{Effect of general viscosity in fluid dynamics}
 
Following our earlier work \cite{JLfive} we now return to the general diffusion approximation for strictly 
hyperbolic systems with general diffusion matrix $B(u)$ of the form
\be
\label{system}
-y u' + A(u) u' =\eps (B(u)u')', 
\ee
and we consider an associated generalized eigenvalue problem.
 
Let $\mu_j$, $\hatl_j$, $\hatr_j$, $j=1,2,\ldots,N$ be 
the eigenvalues and left and right-eigenvectors given by 
\be
\label{eigen}
\aligned
& \bigl(- y + A(v) \bigr) \, \hatr_j(v,y) = \mu_j (v,y) \,  B(v) \, \hatr_j(v,y), 
\\
& \hatl_j(v,y) \cdot \bigl( - y + A(v) \bigr) = \mu_j(v,y) \, \hatl_j(v,y)
\, B(v). 
\endaligned 
\ee
Let us impose the normalization
$$
\hatl_i(v,y) \, B(v) \, \hatr_j(v,y) = 0  \, \text{ if } i \neq j;  
\quad 
\hatl_j(v,y) \, B(v) \, \hatr_j(v,y) = 1.
$$  

In the special case $B(u) = I$ we have simply 
\be
\label{special}
\mu_j(v,y)         = - y + \lam_j(v), 
\quad \hatr_j(v,y) = r_j(v), 
\quad \hatl_j(v,y) = l_j(v).
\ee 
So, we expect the eigenvalues and eigenvectors in \eqref{eigen} to be close to 
to that of \eqref{special}, at least in the case that  
$B(v) = I + \eta \, T(v,\eta)$ where $T(v, \eta) = \bigl( t_{ij} (v, \eta) \bigr)$ 
with $t_{ij}(v, \eta) = O(1)$ and $\eta<<1$. 

\begin{lemma} If $|B(v) - I| = O(\eta)$ is sufficiently small,
then 
\be
\label{approx}
\aligned
& \mu_j(v,y) = - y + \lambda_j(u) + O(\eta),
  \quad \hatl_j(v,y) = l_j(v) + O(\eta),
  \\
& \hatr_j(v,y) = r_j(v) + O(\eta),
 \quad 
 \hatl_i(v,y) \,B(v) \, \del_y \hatr_j(v,y) = O(\eta).
\endaligned
\ee
\end{lemma}

\begin{proof}
When $B$ has the form $B(v) = I + \eta \, T(v)$ where the matrix $T(v)$ 
does not depend on $\eta$, the desired estimates follow from classical results. 
In the slightly more general case when $T= T(v, \eta)$ we may argue as follows. 
First, note that the eigenvalues under consideration are given by 
\be
\label{det}
\det \Big( \bigl(- y \, I + A(v) \bigr) - \mu \, B(v) \Big) = 0.
\ee
Clearly, \eqref{det} is a polynomial equation in $\mu$ whose coefficients
are polynomials in the $t_{ij}$'s. The leading coefficient $det \Bigl( I + \eta \, \bigl( t_{ij} \bigr) \Big)$
is different from zero. Since the roots of \eqref{det} are distinct for
$t_{ij} = 0$, then for $\eta$ sufficiently small the roots of \eqref{det} 
(namely the $\mu_j$'s) are distinct and depend smoothly upon the $t_{ij}$'s with, in addition, 
\be
\label{expan61}
\mu_j(v,y) = - y + \lam_j(v) + \eta \sum_i \lam_{ij}(v,y) \, t_{ij}(v,\eta) + \ldots, 
\ee
where the coefficients $\lam_{ij}$ are smooth functions in $y$.

Since the corresponding left- and right-eigenvectors (which need not 
be normalized at this stage) $\hatl_j$ and $\hatr_j$ are polynomials in 
$a_{ij}$ and $\mu_j$, it follows that they are also smooth in the $a_{ij}$'s. 
Substituting them in \eqref{special} we get
\be
\label{expan62}
\aligned
&\hatl_j(v,y) = l_j(v) + \eta \sum_i l_{ij}(v,y) \, t_{ij}(v,\eta) +
\ldots,\\ 
&\hatr_j(v,y) = r_j(v) + \eta \sum_i r_{ij}(v,y) \, t_{ij}(v,\eta) +
\ldots,
\endaligned
\ee
where $l_{ij}$ and $r_{ij}$ depend smoothly upon $y$. From \eqref{expan62}, it follows
that for $\eta$ sufficiently small the expansion holds for the normalized 
vectors $\hatl_j(v,y)$ and $\hatr_j(v,y)$. Finally, the required estimates \eqref{approx} follow 
from \eqref{expan61} and \eqref{expan62}. This completes the proof of \eqref{approx}. 
\end{proof} 

By the same technique as above we can also prove: 

\begin{theorem} 
Under the conditions that $|D_u B(u)| < \eta$ and the characteristic fields associated with $A(u)$ 
are genuinely nonlinear, the generalized Lax shock inequalities
\[
\hat\lambda_j(u(y+),y)\leq y \leq \hat\lambda_j(u(y-),y) 
\]
are equivalent to the standard Lax shock inequalities, 
\[
\lambda_j(u(y+))\leq y \leq \lambda_j(u(y-)), 
\]
where  $\hatlam_j = \la \hatr_j, A \, \hatr_j \ra$.
\end{theorem}

\begin{proof} We assume that $\nabla \lambda_j \cdot r_j >0$. 
The key observation is that $\nabla \hatlam_j \cdot \hatr_j$ is also positive 
and therefore the same part of the Hugoniot curve is selected by 
the standard and by the generalized Lax shock inequalities. This is so because 
$\nabla \hatlam_j \cdot \hatr_j = \nabla \lam_j \cdot r_j + O(\eta)$ as we now check. 

When $|D_u B(u)| < \eta$ we see that $D_v T(v,\eta)| \leq C$ and that the expansion \eqref{expan62} is valid
for the $v$-derivatives as well. This completes the proof. 
\end{proof}

We end this section with a further discussion of systems of two equations, and provide some explicit calculations 
which lead to sufficient conditions on the diffusion matrix $B$ allowing us to apply the techniques introduced
earlier in this paper. 
Note first that \eqref{det} can be written as
\be
\label{detg}
\det \Big( \bigl( - y \, I + \diag\bigl( \lam_1(v), \lam_2(v) \bigr) - \mu 
\, L(v) \, B(v) \, R(v) \Big) = 0,
\ee
where $L(u)$ and $R(u)$ are matrices of left- and right-eigenvectors
associated with $A(u)$.
Namely, $L(v):= \bigl( l_1(v), \ldots, l_N(v) \bigr)$
and $R(v) := \bigl( r_1(v), \ldots, r_N(v) \bigr)$. 
Setting $L(v) \, B (v) \, R(v)  = \bigl(b_{ij}(v) \bigr)$,
\eqref{detg} becomes 
$$
\bigl( \det B(v) \bigr) \, \mu^2 -  \Big( b_{11}(v) \, (\lam_2 - y) 
+ b_{22}(v) \, (\lam_1 - y) \Big) \, \mu
+ (\lam_1 - y) \, (\lam_2 - y) = 0. 
$$
Solving this quadratic equation in $\mu$ and using the notation
$$
\beta(v) := {b_{12}(v) \, b_{21}(v) \over b_{11}(v) \, b_{22}(v)},
\quad 
a_i(v,y) := {\lam_i(v) - y \over b_{ii}(v)},
$$
we arrive at
$$
\aligned
\mu_1(v,y)
& = {\bigl( a_2(v,y) + a_1(v,y) \bigr) - \bigl( a_2(v,y) - a_1(v,y) \bigr) \,
\Big( 1 + 4 \, \beta(v) \, {a_2(v,y) \, a_1(v,y) \over a_2(v,y) - a_1(v,y)  } \Big)^{1/2} 
\over 2 \, (1 - \beta(v)}, 
\\ 
\mu_2(v,y)
& = {\bigl( a_2(v,y) + a_1(v,y) \bigr) + \bigl( a_2(v,y) - a_1(v,y) \bigr) \,
\Big( 1 + 4 \, \beta(v) \, {a_2(v,y) \, a_1(v,y) \over a_2(v,y) - a_1(v,y) 
}\Big)^{1/2} 
\over 2 \, (1 - \beta(v)}.
\endaligned
$$
It is easy to see that $\mu_1(v)$ and $\mu_2(v)$ are real if 
$0 < \beta(v) < 1$, which is guaranteed if 
\be
\label{ineq}
\aligned 
& b_{11}(v) > 0, \quad b_{22}(v) > 0,
\\
& b_{12} (v) \, b_{21} (v) > 0. 
\endaligned  
\ee
The first condition above is also a necessary condition for viscous shocks 
to be strictly stable in the sense of Majda and Pego \cite{MajdaPego}.

Assuming \eqref{ineq}, an easy calculation shows that the eigenvalues
$\mu_i(v)$ are separated: 
$$
\mu_2(v) - \mu_1(v) \geq C > 0.
$$
It is easy to compute the corresponding right- and left-eigenvectors, namely 
$$
\aligned
& \hatr_i(v,y) = \Big( - b_{12}(v) \, \mu_i(v,y), \lam_1(v) -y - b_{11}(v) 
\, \mu_i(v,y) \Big)^t, 
\\
& \hatl_i(v,y) = \Big( - b_{21}(v) \, \mu_i(v,y), \lam_1(v) -y - b_{11}(v) 
\, \mu_i(v,y) \Big).
\endaligned 
$$
Finally, it is easily seen that the estimates \eqref{approx} hold true when $|B(u) - I| = O(\eta)$.

\section*{Acknowledgments} 
 
KTJ and PGL were supported by a grant (Number 2601-2) from the Indo-French Centre 
for the Promotion of Advanced Research, IFCPAR (Centre Franco-Indien pour la Promotion de 
la Recherche Avanc\'ee, CEFIPRA). 
PGL was also supported by Centre National de la Recherche Scientifique (CNRS) and 
by the A.N.R. grant 06-2-134423 entitled ``Mathematical Methods in General Relativity'' (MATH-GR). 

%************************************************************************************************
 
 \newcommand{\auth}{\textsc}

\end{document}